\documentclass[12pt]{amsart}
\usepackage{amssymb,amsmath, amsthm, enumerate, hyperref}
\usepackage[mathscr]{euscript}
\usepackage[a4paper]{geometry}
\usepackage[T1]{fontenc} 
\usepackage[utf8]{inputenc}   

\theoremstyle{plain}
\newtheorem{theorem}{Theorem}[section]
\theoremstyle{definition}

\numberwithin{equation}{section}

\def\R {\mathbb{R}}

\def\C {\mathbb{C}}

\usepackage{enumitem}
	\setenumerate{leftmargin=0.42cm}

\usepackage{mathtools}
\DeclarePairedDelimiter\floor{\lfloor}{\rfloor}

\def \au {\rm}
\def \ti {\it}
\def \jou {\rm}
\def \bk {\it}
\def \no#1#2#3 {{\bf #1} (#3), #2.}
\def \eds#1#2#3 {#1, #2, #3.}
\newcommand{\T}{(T(t))_{t\ge 0}}
\newcommand{\Tris}{(T_\alpha(t))_{t\ge 0}}
\DeclareMathOperator*{\re}{Re}

\begin{document}

\title[Growth rates under resolvent bounds]
{A remark on the growth rate of operator semigroups under resolvent bounds}

\author[F. Dell'Oro]
{Filippo Dell'Oro}
\address{Politecnico di Milano -- Dipartimento di Matematica
\newline\indent
Via Bonardi 9, 20133 Milano, Italy}
\email{filippo.delloro@polimi.it}

\subjclass[2020]{47D06, 47A10}
\keywords{$C_0$-semigroup, growth rate, resolvent bound, Kreiss condition}

\begin{abstract}
We provide a growth bound for the operator norm of $C_0$-semigroups on Hilbert
spaces under a corresponding growth bound on the resolvent of the semigroup generator.
For some super-linear resolvent growths,
our estimate is sharper than the ones currently available in the literature.
\end{abstract}

\maketitle

\section{Introduction and Main Result}

\noindent 
Let $A$ be the infinitesimal generator of a $C_0$-semigroup $\T$ on a complex Hilbert space $H$.
Denote by $\rho(A)$ the resolvent set of $A$ and by $R(\lambda,A)=(\lambda- A)^{-1}$ the 
resolvent operator at $\lambda\in\rho(A)$.   
We say that $A$ satisfies the {Kreiss condition} if the open right half-plane 
$\C^{+} = \{ \lambda \in \C : \re \lambda > 0\}$
is contained in $\rho(A)$ and there exists a constant $C\geq 1$ such that
\begin{equation}
\label{kreiss}
\|R(\lambda,A)\|\leq \frac{C}{\re \lambda},\quad\,\lambda \in \C^{+}.
\end{equation}
The validity of the Kreiss condition with $C=1$
is equivalent to $\T$ being a contraction semigroup, while for $C>1$ it implies
that $\|T (t)\|=O(t)$ as $t\to\infty$; see \cite{EZ}. Moreover, for every fixed $\gamma \in (0,1)$, there exists
a semigroup on a Hilbert space whose norm grows at least as $t^\gamma$ for $t\to\infty$ and whose infinitesimal
generator satisfies the Kreiss condition; see again \cite{EZ}.
The question whether such examples exist even for $\gamma=1$ has been open for some 
time, as explicitly stated in \cite{EZ,RO}, and was solved negatively in \cite{ARN}. In that article, 
by adapting a technique devised in \cite{CCEL} for discrete semigroups, 
the author proved that \eqref{kreiss} implies the
sharper bound
\begin{equation}
\label{arnold-est}
\|T (t)\|=O\left(\frac{t}{\sqrt{\log t}}\right), \quad\, t\to\infty.
\end{equation}
More details on the Kreiss condition and its applications may be found in \cite{ARN2,EZ2,ROVE,ROVE2}.
We further mention the original article of Kreiss \cite{kreiss-paper}
which contains the first explicit statement of \eqref{kreiss}
for $A$ being a matrix.

In this paper, we consider more general resolvent bounds of the form
$$
\|R(\lambda,A)\|\leq g \Big(\frac{1}{\re \lambda}\Big),\quad\,\lambda \in \C^{+},
$$
where $g:(0,\infty)\to (0,\infty)$ is a non-decreasing function with $g(t)\to\infty$
for $t\to\infty$. 
As in the case of the Kreiss condition, one expects that the growth rate of $g$
reflects on the growth rate of $\|T(t)\|$, and this insight has been 
confirmed in \cite{Bou,HEL,HEL2,HEL3,ROVE}.
In particular, the resolvent bound above implies that 
\begin{equation}
\label{rove-est}
\|T (t)\|=O(g(t)), \quad\, t\to\infty;
\end{equation}
see \cite[Theorem 3.4]{ROVE}. The aim of this short article is to prove the following result.

\begin{theorem}
\label{KREISSLOG}
Let $A$ be the infinitesimal generator of a $C_0$-semigroup $\T$ on a complex Hilbert space $H$
satisfying $\,\C^{+} \subset \rho(A)$.
Assume that there exists a function $g:(0,\infty)\to (0,\infty)$ such that
\begin{equation}
\label{mia-kreiss}
\|R(\lambda,A)\|\leq g \Big(\frac{1}{\re \lambda}\Big),\quad\,\lambda \in \C^{+}.
\end{equation}
Moreover, assume that the function $h(t)=g(t)/t$ is non-decreasing. Then
\begin{equation}
\label{mia-est}
\|T (t)\|= O\left(g(t)\frac{h(t)}{\sqrt{\log t}}\right), \quad\, t\to\infty.
\end{equation}
\end{theorem}

Note that \eqref{mia-kreiss} implies the Kreiss condition \eqref{kreiss}
whenever $h(t)$ is bounded, and in such a case
\eqref{mia-est}
boils down to \eqref{arnold-est}.
Instead, if $h(t)\to\infty$ as $t\to\infty$, then \eqref{mia-est} is sharper than \eqref{rove-est} whenever $h(t)=o(\sqrt{\log t})$ as $t\to\infty$. 
The following two examples illustrate this possibility.
 
\begin{enumerate}
\item[$\bullet$] 
Let $g(t) = t\hspace{0.3mm} \log^\beta (t+1)$ with $\beta \in (0,\frac12)$. In this situation,
Theorem \ref{KREISSLOG} yields the improved estimate 
$$\|T(t)\| = O(t\hspace{0.3mm} \log^{\frac{4\beta-1}{2}} t), \quad\, t\to\infty.$$

\item[$\bullet$]
Let $g(t) = t\hspace{0.3mm} \log\log (t+e)$. In this situation,
Theorem \ref{KREISSLOG} yields the improved estimate 
$$\|T(t)\| = O\left(\frac{t\hspace{0.3mm}\log\log t}{\sqrt{\log t}}\right), \quad\, t\to\infty.$$
\end{enumerate}

The interested reader will have no difficulty in constructing other examples, for instance by
considering iterated logarithms. 

\section{Proof of Theorem \ref{KREISSLOG}}
\label{sec:proof}

\noindent
We begin by noting that since $h(t)$ is non-decreasing then $g(t)$ is non-decreasing and such that $g(t)\to\infty$
for $t\to\infty$. 
Let $s(A) = \sup \big\{{\rm Re}\, \lambda : \lambda \in \sigma(A)\big\}$ be the spectral bound of $A$,
where $\sigma(A)$ denotes the spectrum of $A$. Let also
$$
s_0(A)= \inf \big\{r> s(A) :\, \exists\hspace{0.3mm} C_r>0 : 
\,\|R(\lambda,A)\|\leq C_r\,\text{ whenever } \re \lambda>r\big\}
$$
be the pseudo-spectral bound of $A$. Being $g(t)$ non-decreasing, 
assumption \eqref{mia-kreiss} ensures that $s_0(A)\leq0$.
Hence, calling
$$
\omega_0 = \inf \big\{ \omega \in \R :\, \exists\hspace{0.3mm} M_\omega\geq1 : 
\,\|T(t)\| \leq M_\omega\hspace{0.3mm} {\rm e}^{\omega t}
\,\text{ for all } t\geq0 \big\}
$$
the growth bound of $\T$, the Gerhart-Pr\"uss theorem yields $\omega_0 = s_0(A)\leq0$; 
see e.g.\ \cite[Theorem 5.2.1]{BATTY}.
Therefore, the following integral representation of the resolvent holds
\begin{equation}
\label{int-res}
R(\lambda,A)x=\int_0^\infty {\rm e}^{-\lambda t} T(t)x\,\mathrm{d} t,\quad\, x\in H;\,\ \re\lambda>0.
\end{equation}
For $\alpha>0$, we denote by $\Tris$  the rescaled semigroup given by 
$T_\alpha(t) = {\rm e}^{-\alpha t} T(t)$ when $t\ge0$, and we shall extend $\Tris$ by zero to $\R$. 
Equality~\eqref{int-res} tells that the function $s\mapsto R(\alpha + i s,A)x$ 
is nothing but the Fourier transform of the function
$t\mapsto T_\alpha(t) x$.  
Therefore, applying Plancherel's theorem for square-integrable
functions taking values in a Hilbert space, we deduce that
\begin{equation}
\label{PLANCH}
\int_0^{\infty} \|T_\alpha(\tau)x\|^2\, \mathrm{d}\tau 
= \frac{1}{2\pi}\int_{-\infty}^{\infty} \|R(\alpha + i s,A)x \|^2\, \mathrm{d}s;
\end{equation}
see e.g.\ \cite[Appendix C]{ENG}.
We stress that the validity of \eqref{int-res} for $\re\lambda>0$, as well as the validity of
\eqref{PLANCH} for all $\alpha>0$, is a consequence of $\omega_0\leq 0$. 

\smallskip
As shown in \cite{EZ}, for every $r \in (0,1)$ and every $x,y\in H$, there exists a constant
$K>0$ independent of $r,x,y$ such that
\begin{align}
\label{ez-stima1}
\int_{-\infty}^{\infty}\|R(r + i s,A)x\|^2\, \mathrm{d}s\leq K^2\Big(1+g\Big(\frac1r\Big)\Big)^2\|x\|^2\\
\noalign{\vskip1mm}
\label{ez-stima2}
\int_{-\infty}^{\infty}\|R(r + i s,A^*)y\|^2\, \mathrm{d}s\leq K^2\Big(1+g\Big(\frac1r\Big)\Big)^2\|y\|^2
\end{align}
where $A^*$ is the adjoint of $A$. Actually, these estimates have been obtained in \cite{EZ} only 
for $g(t)=C\hspace{0.2mm} t$ with $C\geq1$, but the same argument applies 
to any non-decreasing function; see also \cite{Bou,EZ2}.
Nevertheless, in order to make the
paper self-contained, we report here the full proof. Let
$\omega,M_\omega\geq1$ 
be such that $\|T(t)\|\leq M_\omega\hspace{0.3mm} {\rm e}^{t(\omega-1)}$
for $t\geq0$. From the resolvent identity and \eqref{mia-kreiss} it follows
that, for $r \in (0,1)$, $s\in\R$ and $x\in H$
\begin{align*}
\|R(r + i s,A)x\|&=\|R(\omega + i s,A)x + (\omega-r)R(r + i s,A)R(\omega + i s,A)x\|\\
&\leq \omega\Big(1+ g\Big(\frac1r\Big)\Big)\|R(\omega + i s,A)x\|.
\end{align*}
Moreover, exploiting \eqref{PLANCH} with $\alpha=\omega$, we deduce that
$$
\int_{-\infty}^{\infty} \|R(\omega + i s,A)x \|^2\, \mathrm{d}s=2 \pi
\int_0^{\infty} \|T_\omega(\tau)x\|^2\, \mathrm{d}\tau \leq \pi M_\omega^2\hspace{0.3mm}\|x\|^2.
$$
Combining the two estimates above, we reach \eqref{ez-stima1} with 
$K=\sqrt{\pi}\hspace{0.3mm} \omega M_\omega$. In order to prove \eqref{ez-stima2}, we
apply the same argument 
to $A^*$ and the adjoint semigroup $(T^*(t))_{t\ge0}$,
recalling that $\|R(\lambda,A^*)\|=\|R(\lambda,A)\|$
for $\lambda\in\rho(A)$ and $\|T(t)\|=\|T(t)^*\|$ for $t\geq0$. 

\smallskip
Next, we show the following bounds
\begin{align}
\label{stima1}
\int_{0}^{t}\|T(\tau)x\|^2\, \mathrm{d}\tau \leq \frac{{\rm e}^2 K^2}{2 \pi}\big(1+g(t)\big)^2\|x\|^2\\
\noalign{\vskip1mm}
\label{stima2}
\int_{0}^{t}\|T^*(\tau)y\|^2\, \mathrm{d}\tau\leq \frac{{\rm e}^2 K^2}{2 \pi}\big(1+g(t)\big)^2\|y\|^2
\end{align}
for every $t>1$ and every $x,y\in H$.
To this end, making use \eqref{PLANCH} with  $\alpha=r\in(0,1)$, we have
$$
{\rm e}^{-2rt} \int_{0}^{t}\|T(\tau)x\|^2\, \mathrm{d}\tau \leq
\int_0^{\infty} \|T_r(\tau)x\|^2\, \mathrm{d}\tau =
\frac{1}{2\pi}\int_{-\infty}^{\infty} \|R(r + i s,A)x \|^2\, \mathrm{d}s.
$$
Owing to \eqref{ez-stima1}, the latter yields
$$
\int_{0}^{t}\|T(\tau)x\|^2\, \mathrm{d}\tau \leq \frac{{\rm e}^{2rt}K^2}{2\pi}
\Big(1+g\Big(\frac1r\Big)\Big)^2\|x\|^2,
$$
and choosing $r = 1/t$ we arrive at \eqref{stima1}. 
In the light of \eqref{ez-stima2}, the same argument applied to 
$(T^*(t))_{t\ge0}$ leads to \eqref{stima2}.

\smallskip
We have now all the ingredients to finish the proof. To this end,
let $t\geq2$ and $x,y\in H$ be arbitrarily fixed.
For every $1\leq a<b\leq t$, we write
\begin{align*}
(b-a)|\langle T(t) x, y \rangle| &= \int_a^b |\langle T(t-\tau) x, T^*(\tau) y \rangle|\mathrm{d}\tau 
\\
&\leq 
\Big(\int_a^b \| T^*(\tau) y\|^2 \mathrm{d}\tau  \Big)^{1/2}
\Big(\int_a^b \| T(t-\tau) x\|^2 \mathrm{d}\tau  \Big)^{1/2}\\
&\leq \frac{{\rm e} K}{\sqrt{2 \pi}}\big(1+g(b)\big)
\bigg(\int_{t-b}^{t-a} \| T(\tau) x\|^2 \mathrm{d}\tau\bigg)^{1/2}\|y\|,
\end{align*}
where the second inequality follows from \eqref{stima2}. Therefore
\begin{equation}
\label{equaz}
\Big(\frac{b-a}{b}\Big)|\langle T(t) x, y \rangle|
\leq \frac{{\rm e} K}{\sqrt{2 \pi}}\big(1+h(b)\big)
\Big(\int_{t-b}^{t-a} \| T(\tau) x\|^2 \mathrm{d}\tau\Big)^{1/2}\|y\|.
\end{equation}
Taking the supremum over $\|y\|=1$, and since $h(b)\leq h(t)$, we get
$$
\Big(\frac{b-a}{b}\Big)\|T(t) x\|
\leq \frac{{\rm e} K}{\sqrt{2 \pi}} \big(1+h(t)\big)
\Big(\int_{t-b}^{t-a} \| T(\tau) x\|^2 \mathrm{d}\tau\Big)^{1/2}.
$$
Set now $N= \floor*{\log t / \log 2}$ where $\floor*{\,\cdot\,}$ is the floor function.
Choosing $a= 2^{n-1}$ and $b=2^{n}$ for $1\leq n\leq N$, we find
$$
\|T(t) x\|^2 \leq \frac{2 {\rm e}^2 K^2}{\pi}\big(1+h(t)\big)^2
\int_{t-2^{n}}^{t-2^{n-1}} \| T(\tau) x\|^2 \mathrm{d}\tau.
$$
As a consequence
\begin{align*}
N\|T(t) x\|^2  
&\leq \frac{2 {\rm e}^2 K^2}{\pi}\big(1+h(t)\big)^2
\sum_{n=1}^{N}\int_{t-2^{n}}^{t-2^{n-1}} \| T(\tau) x\|^2 \mathrm{d}\tau\\
&\leq \frac{2 {\rm e}^2 K^2}{\pi}\big(1+h(t)\big)^2
\int_{0}^{t} \| T(\tau) x\|^2 \mathrm{d}\tau.
\end{align*}
In the light of \eqref{stima1}, we end up with
\begin{equation}
\label{final}
\|T(t) x\| \leq \frac{{\rm e}^2 K^2}{\pi}\frac{1}{\sqrt{\floor*{\log t / \log 2}}}\big(1+h(t)\big)\big(1+g(t)\big)
\|x\|.
\end{equation}
Taking the supremum over $\|x\|=1$, we reach the desired estimate \eqref{mia-est}.
\qed

\section{Concluding Remarks}

\noindent
$\bullet$
In \cite[Theorem 3.4]{ROVE} the authors keep track of the implicit constant appearing in
estimate \eqref{rove-est} by means of the constants $\omega$ and $ M_\omega$ 
used in the proof of Theorem \ref{KREISSLOG}.
The same can be done for the implicit constant appearing in \eqref{mia-est}, 
simply recalling that the constant $K$ in \eqref{final} 
is equal to $\sqrt{\pi}\hspace{0.3mm} \omega M_\omega$.

\medskip
\noindent
$\bullet$
The proof of Theorem \ref{KREISSLOG} can be adapted to give
another proof of estimate~\eqref{rove-est} from \cite{ROVE}. Indeed,
choosing $b=2$ and $a=1$ in \eqref{equaz}, we have for $t\geq2$
\begin{align*}
|\langle T(t) x, y \rangle|
&\leq \frac{\sqrt{2}{\rm e} K}{\sqrt{\pi}}\big(1+h(2)\big)
\Big(\int_{t-2}^{t-1} \| T(\tau) x\|^2 \mathrm{d}\tau\Big)^{1/2}\|y\|\\
&\leq \frac{\sqrt{2}{\rm e} K}{\sqrt{\pi}}\big(1+h(2)\big)
\Big(\int_{0}^{t} \| T(\tau) x\|^2 \mathrm{d}\tau\Big)^{1/2}\|y\|.
\end{align*}
Owing to \eqref{stima1}, we find
$$
|\langle T(t) x, y \rangle| \leq \frac{{\rm e}^2 K^2}{{\pi}}\big(1+h(2)\big)\big(1+g(t)\big)
\|x\|\|y\|,
$$
and taking the supremum over $\|x\|=\|y\|=1$ we arrive at \eqref{rove-est}.
Note that here the monotonicity of the function $h(t)$ is not needed. 
However, it is worth remarking that
this argument is valid only on Hilbert spaces,
while in \cite{ROVE} estimate \eqref{rove-est} has been achieved even for certain classes
of semigroups acting on Banach spaces, such as asymptotically analytic semigroups or
positive semigroups on $L^p$-spaces with $p\neq2$.  


\end{document}